\documentclass[12pt]{amsart}

\usepackage[margin=1in]{geometry}          
\usepackage{amsmath, amscd, amssymb}     
\usepackage{amsfonts}    
\usepackage{amsthm}
\usepackage[all]{xy}               

\def\m{\medskip}

\newtheorem*{theorem}{Theorem}

\newcommand\theoref{Theorem~\ref}

\def\Z{\mathbb Z}
\def\cat{\operatorname{cat}}
\def\TC{\operatorname{TC}}

\def\ts{\times}
\def\p{{\bf Proof. }}

\long\def\forget#1\forgotten{} % 

\begin{document}

\title[On topological complexity of  Eilenberg-MacLane spaces]{On topological complexity of Eilenberg-MacLane spaces }

\author{Yuli Rudyak} 
\address{Department of Mathematics, 
University of Florida,
Gainesville, Florida 32608}
\email{rudyak@ufl.edu}  
 
 
\begin{abstract}
 We note that, for any natural $k$ and every natural $l$ between $k$ and $2k$, there exists a group $\pi$ with $\cat K(\pi,1)=k$ and $\TC(K(\pi,1))=l$. Because of this, we can set up a problem of searching of purely group-theoretical description of $\TC(K(\pi,1))$ as an invariant of $\pi$.
\end{abstract}

\maketitle

\m Below $\cat X$ denotes the Lusternik-Schnirelmann category (normalized, i.e $\cat S^n=1$, see~\cite{CLOT}). Furthermore, we denote by $\TC(X)$ the topological complexity of $X$ defined by Farber~\cite{F1}, but we use the normalized version as ~\cite{R}. 

\m Because of  results of Dranishnikov~\cite[Lemma 2.7 and Theorems 3.6]{D}, we get the following inequalities:
\begin{equation}\label{eq:main}
\cat(G\ts H)\le \TC(G\vee H)\le \cat G+\cat H
\end{equation}

\m Farber asked about calculation of $\TC(K(\pi,1)$'s. It is known that $\cat X \le \TC(X)\le \cat (X\ts X)$ for all $X$,~\cite{F1}. The following observation tells us that, in the class of $(K(\pi,1)$-spaces, the above mentioned inequality gets no new bounds.

 \begin{theorem}\label{t:main} For every natural $k$ and every natural $l$ with $l\le k\le 2k$ there exists a discrete group $\pi$ such that $\pi$ with $\cat K(\pi,1)=k$ and $\TC(K(\pi,1))=l$.  In fact, we can put $\pi=\Z^k*\Z^{l-k}$. 
\end{theorem}

\p Let $T^m$ be the $m$-torus. Then $\cat T^m=m$. Put $r=l-k$ and consider the free product $\pi:=\Z^k*\Z^r$. Then $K(\pi,1)=T^k\vee T^r$, because $\cat (X\vee Y)=\max(\cat X, \cat Y)$ (for good enough spaces X,Y, like CW spaces) . So, $\cat (K(\pi,1))=k$. On the other hand, because of \eqref{eq:main} we have
\[
l=\cat(T^l)=\cat (T^k\ts T^r)\le \TC(T^l\vee T^r) = \TC(K(\pi,1))\le \cat T^k+\cat T^r=k+r=l.
\]
Thus, $\TC(K(\pi,1)=l$.
\qed

\m The TC of groups $\Z^k*\Z^r$ appeared (implicitly) also in~\cite{CP}.

\m Note that the invariant $\cat (K(\pi,1))$ has a known 
purely group-theoretical description. In fact, $\cat K(\pi,1)$ is equal to the cohomological dimension of $\pi$~,\cite{EG}. Now, in view of \theoref{t:main}, the problem of describing of  $\TC (K(\pi,1))$ in purely group-theoretical terms turns out to be essential.

\m {\bf Acknowledgments:} I am grateful to Peter Landweber for his help. The work was partially supported by a grant from the Simons Foundation (\#209424 to Yuli Rudyak).

\end{document}